# NEW INEQUALITIES FOR CONVEX FUNCTIONS
# KONVEKS FONKSİYONLAR İÇİN YENİ EŞİTSİZLİKLER

Mevlüt TUNÇ[1*] ve S. Uğur KIRMACI[2]

[1] *Kilis 7 Aralık Üniversitesi, Fen Edebiyat Fakültesi Matematik Bölümü, Kilis*

[2] *Atatürk Üniversitesi, KKEF, Matematik Eğitimi A.B.D., 25240, Erzurum*

**ABSTRACT**

In the present paper we establish some new integral inequalities analogous to the well known Hadamard's inequality by using a fairly elementary analysis.

**ÖZET**

Bu makalede biz temel analiz işlemlerini kullanarak literatürde iyi bilinen Hadamard eşitsizliğine benzer yeni integral eşitsizlikleri kurduk.

## 1.INTRODUCTION

The following inequality [see Dragomir, (1992)]

$$f\left(\frac{a+b}{2}\right) \leq \frac{1}{b-a}\int_a^b f(x)dx \leq \frac{f(a)+f(b)}{2} \quad (1.1)$$

Which holds for all convex functions $f:[a,b]\to \mathbb{R}$ is known in the literature as Hadamard's inequality. Since its discovery in 1893, Hadamard's inequality [see Hadamard, (1893)] has been proven to be one of the most useful inequalities in mathematical analysis. A number of papers have been written on this inequality providing new proofs, noteworthy extensions, generalizations and numerous applications, see the references cited therein. The main purpose of this paper is to establish some new integral inequalities analogous to that of Hadamard's inequality given in (1.1) involving two convex functions. The analysis used in the proof is elementary.

---

* Sorumlu Yazar: mevluttunc@kilis.edu.tr

## 2. MAIN RESULT

We need the following Lemma proved in [see, Pecaric and Dragomir, 104pp, (1991) ] which deals with the simple characterization of convex functions.

**Lemma A:** The following statements are equivalent to a mapping:

$f:[a,b] \to \mathbb{R}$ ;

i) $f$ is convex on $[a,b]$,

ii) for all $x, y$ in $[a,b]$ the mapping $g:[0,1] \to \mathbb{R}$, defined by

$g(t) = f(tx + (1-t)y)$ is convex on $[0,1]$.

For the proof of this Lemma, see [Pecaric and Dragomir, 104pp, (1991)].

Our main result is given in the following theorem.

**Theorem 2.1:** Let $f:[a,b] \to \mathbb{R}$ be a nonnegative and convex function. Then one has the inequality:

$$\frac{2f(a)}{(b-a)^2}\int_a^b (b-x)f(x)dx + \frac{2f(b)}{(b-a)^2}\int_a^b (x-a)f(x)dx$$
$$\leq \frac{1}{b-a}\int_a^b f^2(x)dx + \frac{M(a,b)}{3} \quad (2.1)$$

where $M(a,b) = f^2(a) + f(a)f(b) + f^2(b)$.

**Proof:** Since $f$ is a convex function on $[a,b]$, then we have that
$f(ta + (1-t)b) \leq tf(a) + (1-t)f(b)$

for all $t \in [0,1]$. Using the elementary inequality $G(a,b) \leq A(a,b),\ 0 \leq a,b$, we can conclude that

$$2\{tf(a)f(ta+(1-t)b)+(1-t)f(b)f(ta+(1-t)b)\}$$
$$\leq f^2(ta+(1-t)b)+t^2f^2(a)$$
$$+2t(1-t)f(a)f(b)+(1-t)^2 f^2(b)$$

By the Lemma A $f(ta+(1-t)b)$ is convex on $[0,1]$, it is integrable on $[0,1]$. Integrating both sides of the above inequality over $t$ on $[0,1]$ we get

$$2f(a)\int_0^1 tf(ta+(1-t)b)dt + 2f(b)\int_0^1 (1-t)f(ta+(1-t)b)dt$$
$$\leq \int_0^1 f^2(ta+(1-t)b)dt + f^2(a)\int_0^1 t^2 dt$$
$$+2f(a)f(b)\int_0^1 t(1-t)dt + f^2(b)\int_0^1 (1-t)^2 dt$$

By substituting $x = ta+(1-t)b$, it is easy to observe that

$$\int_0^1 tf(ta+(1-t)b)dt$$
$$= \frac{1}{b-a}\int_a^b \frac{x-b}{a-b}f(x)dx = \frac{1}{(b-a)^2}\int_a^b (b-x)f(x)dx$$

and

$$\int_0^1 (1-t)f(ta+(1-t)b)dt$$
$$= \frac{1}{b-a}\int_a^b \frac{a-x}{a-b}f(x)dx = \frac{1}{(b-a)^2}\int_a^b (x-a)f(x)dx$$

It can be easily checked that

$$\int_0^1 t^2 dt = \int_0^1 (1-t)^2 dt = \frac{1}{3}, \quad \int_0^1 t(1-t)dt = \frac{1}{6}$$

$$\int_0^1 f^2(ta+(1-t)b)dt = \frac{1}{b-a}\int_a^b f^2(x)dx$$

When above equalities are taken into account, the proof is complete.

**Teorem 2.2:** Let $f:[a,b] \to \mathbb{R}$ be a nonnegative and convex function. Then one has the inequality:

$$\frac{1}{b-a}\int_a^b f(x)dx \leq \frac{1}{2}f\left(\frac{a+b}{2}\right)$$
$$+ \frac{1}{4f\left(\frac{a+b}{2}\right)(b-a)}\int_a^b f^2(x)dx + \frac{N(a,b)}{24f\left(\frac{a+b}{2}\right)} \quad (2.2)$$

where $N(a,b) = f^2(a) + 4f(a)f(b) + f^2(b)$.

**Proof:** Since $f$ is a convex function on $[a,b]$, then we have that

$$f\left(\frac{a+b}{2}\right) = f\left(\frac{ta+(1-t)b}{2} + \frac{(1-t)a+tb}{2}\right)$$
$$\leq \frac{f(ta+(1-t)b) + f((1-t)a+tb)}{2}$$

for all $t \in [0,1]$. Similarly as explained in the proof of inequality (2.1) given above, we obtain

$$f\left(\frac{a+b}{2}\right)f(ta+(1-t)b)+f\left(\frac{a+b}{2}\right)f((1-t)a+tb)$$

$$\leq f^2\left(\frac{a+b}{2}\right)+\frac{f^2(ta+(1-t)b)}{4}$$

$$+\frac{f(ta+(1-t)b)f((1-t)a+tb)}{2}+\frac{f^2((1-t)a+tb)}{4}$$

$$\leq f^2\left(\frac{a+b}{2}\right)+\frac{f^2(ta+(1-t)b)+f^2((1-t)a+tb)}{4}$$

$$+\frac{[tf(a)+(1-t)f(b)][(1-t)f(a)+tf(b)]}{2}$$

$$= f^2\left(\frac{a+b}{2}\right)+\frac{f^2(ta+(1-t)b)+f^2((1-t)a+tb)}{4}$$

$$+\frac{t(1-t)\left[f^2(a)+f^2(b)\right]+\left[t^2+(1-t)^2\right]f(a)f(b)}{2}$$

Likewise as explained in the proof of inequality (2.1) given above, we integrate both sides of the above inequality over $t$ on $[0,1]$, we get

$$f\left(\frac{a+b}{2}\right)\int_0^1 f(ta+(1-t)b)dt + f\left(\frac{a+b}{2}\right)\int_0^1 f((1-t)a+tb)dt$$

$$\leq f^2\left(\frac{a+b}{2}\right)\int_0^1 dt + \frac{1}{4}\int_0^1\{f^2(ta+(1-t)b)+f^2((1-t)a+tb)\}dt$$

$$+\frac{f^2(a)+f^2(b)}{2}\int_0^1 t(1-t)dt + \frac{f(a)f(b)}{2}\int_0^1\left(t^2+(1-t)^2\right)dt$$

By substituting $ta+(1-t)b=x$ and $(1-t)a+tb=y$ it can be easily obtained,

$$\int_0^1 f(ta+(1-t)b)dt = \int_0^1 f((1-t)a+tb)dt = \frac{1}{b-a}\int_a^b f(x)dx$$

$$\int_0^1 t(1-t)dt = \frac{1}{6}, \quad \int_0^1\left(t^2+(1-t)^2\right)dt = \frac{2}{3}$$

From the above equalities it is easily obtained that

$$\frac{f\left(\frac{a+b}{2}\right)}{b-a}\int_0^1 f(x)dx + \frac{f\left(\frac{a+b}{2}\right)}{b-a}\int_0^1 f(x)dx = 2\frac{f\left(\frac{a+b}{2}\right)}{b-a}\int_0^1 f(x)dx$$

$$\leq f^2\left(\frac{a+b}{2}\right) + \frac{1}{2(b-a)}\int_0^1 f^2(x)dx + \frac{f^2(a)+f^2(b)}{12} + \frac{f(a)f(b)}{3}$$

Now dividing both sides of the above inequality by $2f\left(\frac{a+b}{2}\right)$ we get the desired inequality in (2.2). The proof is complete.

**Remark:** If we choose $a=0$ and $b=1$ and the convex functions $f(x)=x$, then it is easy to observe that the inequalities obtained (2.1) and (2.2) are certain in the sense that we hold equalities in (2.1) and (2.2).

**Theorem 2.3:** Let $f:[a,b]\to \mathbb{R}$ be a nonnegative and convex function. Then

$$\frac{2}{(b-a)^2}\int_a^b\int_a^b\int_0^1 f(tx+(1-t)y)(tf(x)+(1-t)f(y))dtdydx$$

$$\leq \frac{1}{(b-a)^2}\int_a^b\int_a^b\int_0^1 f^2(tx+(1-t)y)dtdydx \quad (2.3)$$

$$+\frac{2}{3(b-a)}\int_a^b f^2(x)dx + \frac{\psi(a,b)}{12}$$

where $\psi(a,b) = f^2(a) + 2f(a)f(b) + f^2(b)$.

**Proof:** Since $f$ is a convex function on $[a,b]$, then we have that
$$f(tx+(1-t)y) \leq tf(x)+(1-t)f(y)$$
for all $x,y\in[a,b]$ and $t\in[0,1]$. Using the elementary inequality $G(a,b)\leq A(a,b)$, $0\leq a,b$, we get,

$$2f(tx+(1-t)y)(tf(x)+(1-t)f(y))$$
$$\leq f^2(tx+(1-t)y)+t^2 f^2(x)$$
$$+2t(1-t)f(x)f(y)+(1-t)^2 f^2(y)$$

As explained in the proof of inequality (2.1) given above, we integrate both sides of that above inequality over $t$ on $[0,1]$, and we obtain

$$2\int_0^1 f(tx+(1-t)y)(tf(x)+(1-t)f(y))dt$$
$$\leq \int_0^1 f^2(tx+(1-t)y)dt + f^2(x)\int_0^1 t^2 dt$$
$$+2f(x)f(y)\int_0^1 t(1-t)dt + f^2(y)\int_0^1 (1-t)^2 dt \qquad (2.3.1)$$
$$=\int_0^1 f^2(tx+(1-t)y)dt + \frac{1}{3}f^2(x) + \frac{1}{3}f(x)f(y) + \frac{1}{3}f^2(y)$$

Integrating both sides of that above inequality (2.3.1) over $x$ and $y$ on $[a,b]^2$ we obtain

$$2\int_a^b\int_a^b\int_0^1 f(tx+(1-t)y)(tf(x)+(1-t)f(y))dtdydx$$
$$\leq \int_a^b\int_a^b\int_0^1 f^2(tx+(1-t)y)dtdydx$$
$$+\frac{1}{3}\int_a^b\int_a^b f^2(x)dydx + \frac{1}{3}\int_a^b\int_a^b f(x)f(y)dydx + \frac{1}{3}\int_a^b\int_a^b f^2(y)dydx \qquad (2.3.2)$$
$$\leq \int_a^b\int_a^b\int_0^1 f^2(tx+(1-t)y)dtdydx$$
$$+\frac{2(b-a)}{3}\int_a^b f^2(x)dx + \frac{1}{3}\int_a^b f(x)dx\int_a^b f(y)dy$$

By using the right half of the Hadamard's inequality given in (1.1) on the right side of (2.3.2) we write

$$2\int_a^b\int_a^b\int_0^1 f(tx+(1-t)y)(tf(x)+(1-t)f(y))dtdydx$$

$$\leq \int_a^b\int_a^b\int_0^1 f^2(tx+(1-t)y)dtdydx + \frac{2(b-a)}{3}\int_a^b f^2(x)dx$$

$$+\frac{1}{3}\frac{f(a)+f(b)}{2}\cdot\frac{f(a)+f(b)}{2}(b-a)^2$$

$$=\int_a^b\int_a^b\int_0^1 f^2(tx+(1-t)y)dtdydx + \frac{2(b-a)}{3}\int_a^b f^2(x)dx$$

$$+\frac{1}{12}\left(f^2(a)+2f(a)f(b)+f^2(b)\right)(b-a)^2$$

Now multiplying both sides of the above inequality by $1/(b-a)^2$ we get desired inequality in (2.3). The proof is complete.

**Theorem 2.4:** Let $f:[a,b]\to\mathbb{R}$ be a nonnegative and convex function. Then

$$\frac{2}{b-a}\int_a^b\int_0^1 f\left(tx+(1-t)\frac{a+b}{2}\right)\left(tf(x)+(1-t)f\left(\frac{a+b}{2}\right)\right)dtdx$$

$$\leq \frac{1}{b-a}\int_a^b\int_0^1 f^2\left(tx+(1-t)\frac{a+b}{2}\right)dtdx + \frac{\psi(a,b)}{12}(b-a+2) \quad (2.4)$$

where $\psi(a,b)=f^2(a)+2f(a)f(b)+f^2(b)$.

**Proof:** Since $f$ is a convex function on $[a,b]$, then we have that

$$f\left(tx+(1-t)\frac{a+b}{2}\right)\leq tf(x)+(1-t)f\left(\frac{a+b}{2}\right)$$

for all $x,y\in[a,b]$ and $t\in[0,1]$, we get

$$2f\left(tx+(1-t)\frac{a+b}{2}\right)\left(tf(x)+(1-t)f\left(\frac{a+b}{2}\right)\right)$$

$$\leq f^2\left(tx+(1-t)\frac{a+b}{2}\right)+t^2f^2(x)$$

$$+2t(1-t)f(x)f\left(\frac{a+b}{2}\right)+(1-t)^2 f^2\left(\frac{a+b}{2}\right)$$

As explained in the proof of inequality (2.1) given above, we integrate both sides of that above inequality over $t$ on $[0,1]$, and we obtain

$$2\int_0^1 f\left(tx+(1-t)\frac{a+b}{2}\right)\left(tf(x)+(1-t)f\left(\frac{a+b}{2}\right)\right)dt$$
$$\leq \int_0^1 f^2\left(tx+(1-t)\frac{a+b}{2}\right)dt + f^2(x)\int_0^1 t^2 dt \qquad (2.4.1)$$
$$+ 2f\left(\frac{a+b}{2}\right)f(x)\int_0^1 t(1-t)dt + f^2\left(\frac{a+b}{2}\right)\int_0^1 (1-t)^2 dt$$
$$= \int_0^1 f^2\left(tx+(1-t)\frac{a+b}{2}\right)dt + \frac{1}{3}f^2(x) + \frac{1}{3}f\left(\frac{a+b}{2}\right)f(x) + \frac{1}{3}f^2\left(\frac{a+b}{2}\right)$$

As explained in the proof of inequality (2.3) given above, integrating both sides of that above inequality (2.4.1) over $x$ on $[a,b]$, and using the right half of the Hadamard's inequality given in (1.1) and convexity of $f$ we obtain

$$2\int_a^b\int_0^1 f\left(tx+(1-t)\frac{a+b}{2}\right)\left(tf(x)+(1-t)f\left(\frac{a+b}{2}\right)\right)dtdx$$
$$\leq \int_a^b\int_0^1 f^2\left(tx+(1-t)\frac{a+b}{2}\right)dtdx$$
$$+\frac{1}{3}\int_a^b f^2(x)dx + \frac{1}{3}f\left(\frac{a+b}{2}\right)\int_a^b f(x)dx + \frac{(b-a)}{3}f^2\left(\frac{a+b}{2}\right)$$
$$\leq \int_a^b\int_0^1 f^2\left(tx+(1-t)\frac{a+b}{2}\right)dtdx$$
$$+\frac{1}{3}(b-a)^2\left(\frac{f(a)+f(b)}{2}\right)^2$$
$$+\frac{1}{3}\frac{f(a)+f(b)}{2}(b-a)\frac{f(a)+f(b)}{2} + \frac{(b-a)}{3}\left(\frac{f(a)+f(b)}{2}\right)^2$$
$$= \int_a^b\int_0^1 f^2\left(tx+(1-t)\frac{a+b}{2}\right)dtdx$$
$$+\frac{\psi(a,b)}{12}\left((b-a)^2+(b-a)+(b-a)\right) \qquad (2.4.2)$$

Now dividing both sides of (2.4.2) by $(b-a)$ and rewriting (2.4.2) we get the required inequality in (2.4). The proof is complete.

## 3. APPLICATIONS FOR SPECIAL MEANS

As in [Kırmacı, 2004], we shall consider the means as arbitrary real numbers $a, b$, $a \neq b$. In the resources there includes

$A(a,b) = \dfrac{a+b}{2}, \quad a,b \in \mathbb{R}_+$ (arithmetic mean)

$K(a,b) = \sqrt{\dfrac{a^2+b^2}{2}}, \quad a,b \in \mathbb{R}_+$ (quadratic mean)

$L(a,b) = \dfrac{b-a}{\ln|b|-\ln|a|}, \quad |a| \neq |b|,\ ab \neq 0,\ a,b \in \mathbb{R}_+$ (logarithmic mean)

$G(a,b) = \sqrt{ab}, \quad a,b \in \mathbb{R}_+$ (geometric mean).

Now, using the results of Section 2, we illustrate some applications of special means of real numbers.

**Proposition 3.1:** Let $0 < a < b$. Then one has the inequality,

$$\frac{4A(a,b)}{L(a,b)} \leq \frac{2}{3}\frac{K^2(a,b)}{G^2(a,b)} + \frac{10}{3}$$

**Proof:** The proof is immediate from Theorem 2.1 as applied for $f(x) = \dfrac{1}{x}$ and the details are omitted.

**Proposition 3.2:** Let $0 < a < b$. Then one has the inequality,

$$\frac{1}{L(a,b)} \leq \frac{1}{2A(a,b)} + \frac{A(a,b)}{4G^2(a,b)} + \frac{A(a,b)\left(2K^2(a,b)+4G^2(a,b)\right)}{24G^2(a,b)}$$

**Proof:** The assertion follows from Theorem 2.2 as applied to $f(x) = \dfrac{1}{x}$ and the details are omitted.

**REFERENCES**

Dragomir, S.S. (1992) Two mappings in connection to Hadamard's inequalities, J.Math.Anal.Appl. 167, 49-56.

Hadamard, J. (1893).Etude sur les properties des fonctions entieres et en particulier d'une fonction consideree par Riemann, J.Math.Pures appl. 58, 171-215.

Heing, H.P., Maligranda, L. (1991/92) Chebyshev inequality in function spaces, Real Analysis Exchange 17, 211-47.

Kırmacı, U.S. (2004). Inequalities for differentiable mappings and applications to special means of real numbers and to midpoint formula, Appl. Math.Comput. 147, 137−146.

Maligranda, L., Pecaric, J.E. and Persson, L.E. (1994). On some inequalities of the Gruss-Barnes and Borell type, J.Math.Ana.Appl. 187, 306–323.

Mitrinovic, D.S. Analytic Inequalities, Springer-Verlag, Berlin, New York 1970.

Pachpatte, B.G. (2003). On some inequalities for convex functions, RGMIA Res. Rep. Coll., 6(E)

Pecaric, J.E., Dragomir, S.S. (1991). A generalization of Hadamard's inequality for Isotonic linear functionals, Radovi Matematicki 7, 103-107.

****